\documentclass{amsart}
\usepackage{amsfonts}
\usepackage{amsthm}
\usepackage{geometry}
\usepackage{amsmath}
\usepackage{txfonts}

\usepackage{graphicx}
\usepackage{wrapfig}
\usepackage{ascmac}
\usepackage{bm}
\usepackage{indentfirst}
\usepackage{tikz}
\usepackage{color}
\usepackage{ulem}
\usepackage{comment}

\usetikzlibrary{matrix,calc,arrows,decorations.markings}

\newtheorem{thm}{Theorem}[section]
\newtheorem{rem}[thm]{Remark}
\newtheorem{defi}[thm]{Definition}
\newtheorem{ex}[thm]{Example}

\newtheorem{lem}[thm]{Lemma}
\newtheorem{prop}[thm]{Proposition}
\newtheorem{nota}[thm]{Notation}

\def\NZQ{\Bbb}

\def\RR{{\NZQ R}}
\def\CC{{\NZQ C}}

\def\PP{{\NZQ P}}

\def\SS{{\NZQ S}}
\def\TT{{\NZQ T}}

\def\ml{\mathcal{C}}
\def\ml1{\mathcal{C}^1}
\def\mlb1{\mathcal{C}_{b}^{1}}

\def\frk{\frak}

\def\Phi{{\frk n}}

\def\dim{{\rm dim}}
\def\rank{{\rm rank}}

\def\A{{\mathcal A}}
\def\B{{\mathcal B}}

\newcommand{\R}[0]{\mathbb{R}}
\newcommand{\C}[0]{\mathbb{C}}

\newcommand{\Zt}[0]{\mathcal Z}

\usepackage{color}
\usepackage{ifthen}

\newboolean{usecolor}
\setboolean{usecolor}{true}

\ifthenelse{\boolean{usecolor}}
{
  \definecolor{colore}{cmyk}{0,1,0.6,0}
  \definecolor{coloregen}{cmyk}{0.7,0,1,0}
  \definecolor{coloresimo}{cmyk}{1,0.6,0,0}
}
{
  \definecolor{colore}{cmyk}{0,0,0,1}
  \definecolor{coloregen}{cmyk}{0,0,0,1}
  \definecolor{coloresimo}{cmyk}{0,0,0,1}
}

\title{On the non-very generic intersections in discriminantal arrangements}
\author{Simona Settepanella}
\author{So Yamagata}
\address[1]{Department of Economics and Statistics, Torino University,Italy}
\address[2]{Department of Mathematics, Hokkaido University, Japan.}
\email{simona.settepanella@unito.it}
\email{so.yamagata@math.sci.hokudai.ac.jp}

\thanks{The second author was supported by the Program for Leading Graduate Schools (Hokkaido University Ambitious
Leader's Program) and JSPS Research Fellowship for Young Scientists Grant Number 20J10012.}

\begin{document}

\maketitle

\begin{abstract} 
In 1985 Crapo introduced in \cite{Crapo} a new mathematical object that he called \textit{geometry of circuits}. Four years later, in 1989, Manin and Schechtman defined in \cite{MS} the same object and called it \textit{discriminantal arrangement}, the name by which it is known now a days. Those discriminantal arrangements $\B(n,k,\A^0)$ are builded from an arrangement $\A^0$ of $n$ hyperplanes in general position in a $k$-dimensional space and their combinatorics depends on the arrangement $\A^0$. On this basis, in 1997 Bayer and Brandt (see \cite{BB}) distinguished two different type of arrangements $\A^0$ calling \textit{very generic} the ones for which the intersection lattice of $\B(n,k,\A^0)$ has maximum cardinality and \textit{non-very generic} the others. Results on the combinatorics of $\B(n,k,\A^0)$ in the very generic case already appear in Crapo \cite{Crapo} and in 1997 in Athanasiadis \cite{Atha} while the first known result on non-very generic case is due to Libgober and the first author in 2018. In their paper \cite{LS} they provided a necessary and sufficient condition on $\A^0$ for which the cardinality of rank 2 intersections in $\B(n,k,\A^0)$ is not maximal anymore. In this paper we further develop their result providing a sufficient condition on $\A^0$ for which the cardinality of rank r, $r \geq 2$, intersections in $\B(n,k,\A^0)$ decreases.  
\end{abstract}

\section{Introduction}

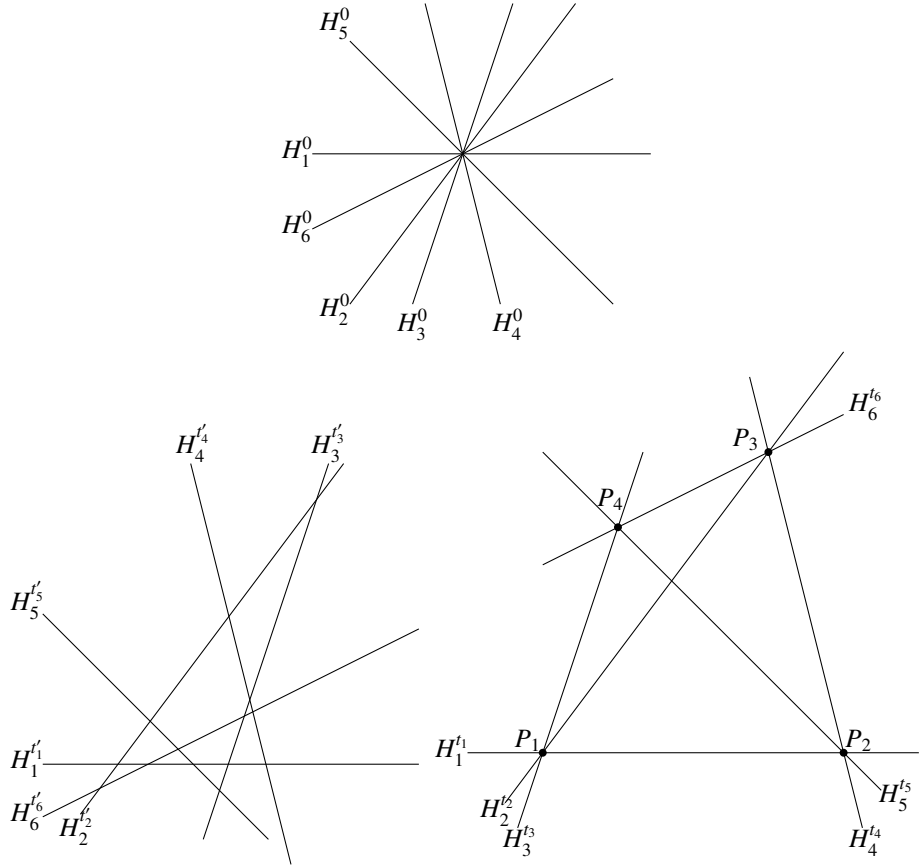
\begin{figure}[h]
\centering
\begin{tikzpicture}
\coordinate (0) at (-2,0);
\coordinate (1) at (-3/2,-2);
\coordinate (2) at (-2/3,-2);
\coordinate (3) at (1/2,-2);
\coordinate (4) at (2,-2);
\coordinate (5) at (5/2,0);
\coordinate (6) at (2,1);
\coordinate (7) at (3/2,2);
\coordinate (8) at (-1/2,2);
\coordinate (9) at (2/3,2);
\coordinate (10) at (-3/2,3/2);
\coordinate (12) at (-2,-1);

\coordinate [label=$H_1^0$] (H1) at (-2.2,-0.3);
\coordinate [label=$H_2^0$] (H2) at (-1.7,-2.4);
\coordinate [label=$H_3^0$] (H3) at (-2/3,-2.6);
\coordinate [label=$H_4^0$] (H4) at (1/2+0.1,-2.6);
\coordinate [label=$H_5^0$] (H5) at (-3/2-0.2,3/2-0.1);
\coordinate [label=$H_6^0$] (H6) at (-2.2,-1.3);

\begin{scope}
\draw (0) -- (5);
\draw (1) -- (7);
\draw (3) -- (8);
\draw (2) -- (9);
\draw (4) -- (10);
\draw (6) -- (12);

\end{scope}
\end{tikzpicture} \\
\begin{tikzpicture}
\coordinate (0) at (-1,0);
\coordinate (1) at (-1/2,-2/3);
\coordinate (2) at (4/3 - 1/5,-1);
\coordinate (3) at (2/3+0.5 - 1/5,4);*
\coordinate (4) at (-1,2);*
\coordinate (5) at (4,0);
\coordinate (6) at (-1,-1/2-0.2);*
\coordinate (7) at (3,4);
\coordinate (8) at (2 - 1/5+0.5,-4/3);*
\coordinate (9) at (3 - 1/5,4);
\coordinate (10) at (2,-1);*
\coordinate (12) at (4,2-0.2);*
\coordinate [label=$H_1^{t'_1}$] (H1) at (-1.2,-0.3);
\coordinate [label=$H_2^{t'_2}$] (H2) at (-0.55,-1.2);
\coordinate [label=$H_3^{t'_3}$] (H3) at (2.8,3.9);
\coordinate [label=$H_4^{t'_4}$] (H4) at (0.5+0.5,3.9);
\coordinate [label=$H_5^{t'_5}$] (H5) at (-1.2,1.3+0.5);
\coordinate [label=$H_6^{t'_6}$] (H6) at (-1.2,-1.0);

\begin{scope}
\draw (0) -- (5); 
\draw (1) -- (7); 
\draw (3) -- (8); 
\draw (2) -- (9); 
\draw (4) -- (10); 
\draw (6) -- (12); 
\end{scope}
\end{tikzpicture}
\begin{tikzpicture}
\coordinate (0) at (-1,0);
\coordinate (1) at (-1/2,-2/3);
\coordinate (2) at (-1/3,-1);
\coordinate (3) at (17/4,-1);
\coordinate (4) at (9/2,-1/2);
\coordinate (5) at (5,0);
\coordinate (6) at (4,9/2);
\coordinate (7) at (4,16/3);
\coordinate (8) at (11/4,5);
\coordinate (9) at (4/3,4);
\coordinate (10) at (0,4);
\coordinate (12) at (0,5/2);
\coordinate (p1) at (0,0);
\coordinate (p2) at (4,0);
\coordinate (p3) at (3,4);
\coordinate (p4) at (1,3);

\coordinate [label=$H_1^{t_1}$] (H1) at (-1.2,-0.3);
\coordinate [label=$H_2^{t_2}$] (H2) at (-0.6,-1.1);
\coordinate [label=$H_3^{t_3}$] (H3) at (-0.3,-1.5);
\coordinate [label=$H_4^{t_4}$] (H4) at (4.3,-1.5);
\coordinate [label=$H_5^{t_5}$] (H5) at (4.7,-0.9);
\coordinate [label=$H_6^{t_6}$] (H6) at (4.3,4.3);
\coordinate [label=$P_1$] (P1) at (-0.2,-0.1);
\coordinate [label=$P_2$] (P2) at (4.2,-0.1);
\coordinate [label=$P_3$] (P3) at (2.7,3.9);
\coordinate [label=$P_4$] (P4) at (0.9,3.1);

\begin{scope}
\draw (0) -- (5);
\draw (1) -- (7);
\draw (3) -- (8);
\draw (2) -- (9);
\draw (4) -- (10);
\draw (6) -- (12);
\fill (p1) circle (1.5pt);
\fill (p2) circle (1.5pt);
\fill (p3) circle (1.5pt);
\fill (p4) circle (1.5pt);

\end{scope}
\end{tikzpicture}
\caption{Central generic arrangement of $6$ lines in $\R^2$, its generic translation on the left and its non-(very) generic translation on the right.}\label{fig:quad}
\end{figure}

In 1989, Manin and Schechtman (\cite{MS}) introduced a family of arrangements of hyperplanes 
generalizing classical braid arrangements, which they called the {\it discriminantal arrangements} (p.209 \cite{MS}). 
Such  an arrangement $\B(n,k,\A^0), n,k \in {\bf N}$ for $k \ge 2$ depends on a choice 
$\A^0=\{H^0_1,...,H^0_n\}$ of a collection of hyperplanes in general position in $\CC^k$, i.e., such that $\dim\bigcap_{i \in K, \mid K\mid=k} H_i^0=0$. It consists of parallel translates of $H_1^{t_1},...,H_n^{t_n}, (t_1,...,t_n) \in \CC^n$ which fail to form a general position arrangement in $\CC^k$. $\B(n,k,\A^0)$ can be viewed as a generalization of the braid arrangement (\cite{OT}) with which $\B(n,1)=\B(n,1,\A^0)$ coincides. \\
These arrangements have several beautiful relations with diverse problems such as the Zamolodchikov equation with its relation to higher category theory (see Kapranov-Voevodsky \cite{KV3}, see also \cite{KV1},\cite{KV2}), the vanishing of cohomology of bundles on toric varieties (\cite{Per}), the representations of higher braid groups (see \cite{Koh}) and, naturally, with combinatorics. The latter is the connection we are mainly interested in and it goes from matroids to special configurations of points, from fiber polytopes to higher Bruhat orders. \\
Manin and Schechtman introduced discriminantal arrangements as higher braid arrangements in order to introduce higher Bruhat orders which model the set of minimal path through a discriminantal arrangement. Even if Ziegler showed (see Theorem 4.1 in \cite{Zie}) in 1991 that we have to choose a cyclic arrangement instead of discriminantal arrangement for this, few years later, in a subsequent work (see \cite{FZ}) Felsner and Ziegler reintroduced the combinatorics of discriminantal arrangement in the study of higher Bruhat orders (this connection uses fiber polytopes as observed by Falk in \cite{Falk}). From a different perspective, unknown in the literature of discriminantal arrangement until Athanasiadis pointed it out in 1999 (see \cite{Atha}), Crapo introduced for the first time in 1985 (see \cite{Crapo}) what he called \textit{geometry of circuits} and which is the matroid $M(n,k, \mathcal{C})$ of circuits of the configuration $\mathcal{C}$ of $n$ generic points in $\RR^k$. The circuits of the matroid $M(n,k, \mathcal{C})$ are the hyperplanes of $\B(n,k,\A^0),$ $\A^0$ arrangement of $n$ hyperplanes in $\RR^k$ orthogonal to the vectors joining the origin with the $n$ points in $\mathcal{C}$ (for further development see \cite{CR}).\\
Both Manin-Schechtman (\cite{MS}) and Crapo (\cite{Crapo}) were mainly interested in the arrangements $\B(n,k,\A^0)$ for which the intersection lattice is constant when $\A^0$ varies within a Zariski open set $\Zt$ in the space of general position arrangements. Crapo shows that, in this case, the matroid $M(n,k)$ is isomorphic to the Dilworth completion of the $k$-th lower truncation of the Boolean algebra of rank $n$. More recently in \cite{Atha}, Athanasiadis proved a conjecture by Bayer and Brandt (see \cite{BB}) providing a full description of combinatorics of $\B(n,k,\A^0)$ when $\A^0$ belongs to $\Zt$. Following \cite{Atha} (more precisely Bayer and Brandt), we call arrangements $\A^0$ in $\Zt$ {\it very generic}, non-very generic otherwise.\\
However Manin and Schechtman do not describe the set $\Zt$ of very generic arrangements explicitly, which, in time, led to the misunderstanding that the combinatorial type of $\B(n,k,\A^0)$ was independent from 
the arrangement $\A^0$ (see for instance, \cite{Orl}, sect. 8, \cite{OT} or \cite{Law}).  Neither Crapo in \cite{Crapo} provided a description of $\Zt$ even if he presented the first known example of a non-very generic arrangement: $6$ lines in generic position in $\RR^2$ which admit translated that are respectively sides and diagonals of a quadrilateral as in Figure \ref{fig:quad} (Crapo calls it a quadrilateral set). Few years later in 1994, Falk provided an higher dimensional example of non-very generic arrangement of $6$ planes in $\RR^3$ (see \cite{Falk}). Similar to Crapo's example, Falk's example too turned out to be related to a special configuration of lines, this time in projective plane (see \cite{SSY1},\cite{SSY2}).\\

\noindent
In 2018 the first general result on non-very generic arrangements is provided. In \cite{LS} Libgober and the first author described a sufficient {\it geometric} condition on the arrangement $\A^0$ to be non-very generic. This condition ensures that $\B(n,k,\A^0)$ admits codimension $2$ strata of multiplicity $3$ which do not exist in very generic case. It is given in terms of the notion of \textit{dependency} for the arrangement $\A_{\infty}$ in $\PP^{k-1}$ of hyperplanes $H_{\infty,1},...H_{\infty,n}$ which are the intersections of projective closures of $H_1^0,...,H_n^0 \in \A^0$ with the hyperplane at infinity. 
Their main result shows that $\B(n,k,\A^0), k>1$ admits a codimension two stratum of multiplicity $3$ if and only if $\A_{\infty}$ is an arrangement in $\PP^{k-1}$ admitting a restriction\footnote{Here restriction is the standard restriction of arrangements to subspaces as defined in \cite{OT}.} which is a dependent arrangement. This construction generalizes Falk's example which corresponds to the case $n=6, k=3$ and which has been object of study in two subsequent papers, (\cite{SSY1}, \cite{SSY2}) by Sawada and the first and second authors. In those papers the authors proved how the arrangement $\A^0$ of $6$ planes in $\RR^3$ (resp. $\CC^3$) for which the rank 2 intersections of $\B(6,3,\A^0)$ are in minimal number corresponds to Pappus's (resp. Hesse's) configuration providing a main example of what conjectured by Crapo that the intersection lattice of the discriminantal arrangement represents a combinatorial way to encode special configurations of points in the space. Notice that in \cite{SSY1} the authors connected the non-very generic arrangements $\A^0$ of $n$ planes in $\CC^3$ to well defined hypersurfaces in Grassmannian $Gr(3,n)$. \\   
In this paper we advance the study of non-very generic arrangements and generalize the dependency condition given in \cite{LS} providing a sufficient condition for the existence in rank $r \geq 2$ of non-very generic intersections, i.e. intersections which doesn't exist in $\B(n,k,\A^0), \A^0 \in \Zt$. In particular we call an intersection of $r$ hyperplanes in $\B(n,k,\A^0)$ which satisfies the following property a \textit{simple} intersection: if the arrangement $\A^0$ is very generic then all simple intersections of multiplicity $r$ have rank $r$ (that is they are $r$ hyperplanes intersecting transversally). Then we provide a geometric necessary and sufficient condition for the existence of simple intersections of multiplicity $r$ in rank strictly lower than $r$, i.e. simple non-very generic intersections. This result firstly connect configurations of non-very generic points to special families of graphs (called $K_{\TT}$-configurations) which help to understand $\B(n,k,\A^0)$ for $\A^0 \notin \Zt$ (as conjectured by Crapo in \cite{Crapo}). Secondly it reduces the geometric problem of the existence of special (non-very generic) configurations of points to a combinatorial problem on the numerical properties that $r$ subsets of indices $L_i \subset \{1, \ldots, n\}, i=1,\ldots , r$ of cardinality $k+1$ have to satisfy in order for the  $K_{\TT}$-configuration, $\TT=\{L_1,\ldots,L_r\}$, to give rise to a simple non-very generic intersection.
The latter problem is left open together with the problem of necessary and sufficient conditions for the existence of intersections in $\B(n,k,\A^0)$ which are nor simple nor very generic. \\
The content of the paper is the following. In Section \ref{sec:prelim}, we recall the definition of discriminantal arrangement and basic properties of the intersection lattice of discriminantal arrangement in very generic case. We also give the definition of \textit{simple} intersection. In Section \ref{sec:geom}, we introduce the notion of $K_\TT$-translated and $K_\TT$-configuration associated to a generic arrangement $\A^0$ providing a geometric condition for $\A^0$ to be non-very generic (Theorem \ref{th:conf}). 

\section{Preliminaries}\label{sec:prelim}
\subsection{Discriminantal arrangement}
Let $H^0_i, i=1,...,n$ be a central arrangement in $\CC^k, k<n$ which is generic\footnote{Notice that, in general, generic, referred to an arrangement of hyperplanes, has a slightly different meaning. With an abuse of notation, we use the word \textit{generic} in this case since the defined property is equivalent to the existence of a translated of the given central arrangement which is generic in the classical sense.}, i.e. any $m$ hyperplanes intersect in codimension $m$ at any point except for the origin for any $m \leq k$. We will call such an arrangement a \textit{central generic arrangement}.
Space of parallel translates $\SS(H_1^0,...,H_n^0)$ (or simply $\SS$ when dependence on $H_i^0$ is clear or not essential) is the space of $n$-tuples of translates $H_1, \dots, H_n$ such that either $H_i \cap H_i^0=\emptyset$ or $H_i=H_i^0$ for any $i=1, \dots, n$.\\
One can identify $\SS$ with $n$-dimensional affine space $\CC^n$ in such a way that $(H^0_1, \dots, H^0_n)$ corresponds to the origin.  In particular, an ordering of hyperplanes in $\A$ determines the coordinate system in $\SS$ (see \cite{LS}). \\
Given a central generic arrangement $\A$ in $\CC^k$ formed by hyperplanes $H_i, i=1, \dots, n$ {\it the trace at infinity}, denoted by $\A_{\infty}$, is the arrangement formed by hyperplanes $H_{\infty,i}= \bar H^0_i \cap H_{\infty}$ in the space $H_{\infty} \simeq \PP^{k-1}(\CC)$, where $\bar H^0_i$ are projective closures of affine hyperplanes $H^0_i$ in the compactification $\PP^k(\CC)$ of $\CC^k\simeq \PP^k(\CC)\setminus H_{\infty}$. Notice that condition of genericity is equivalent to $\bigcup_i  H^0_{\infty,i}$ being a normal crossing divisor in $\PP^{k-1}(\CC)$, i.e. $\A_{\infty}$ is a generic arrangement.\\
The trace $\A_{\infty}$ of an arrangement $\A$ determines the space of parallel translates $\SS$ (as a subspace in the space of $n$-tuples of hyperplanes in $\PP^k$). 
Fixed a generic central arrangement $\A$, consider the closed subset of $\SS$ formed by those collections which fail to form a generic arrangement. This subset of $\SS$ is a union of hyperplanes $D_L \subset \SS$ (see \cite{MS}). Each hyperplane $D_L$ corresponds to a subset $L = \{ i_1, \dots, i_{k+1} \} \subset$  [$n$] $\coloneqq \{ 1, \dots, n \}$ and it consists of $n$-tuples of translates of hyperplanes $H_1^0, \dots, H_n^0$ in which translates of $H_{i_1}^0, \dots, H_{i_{k+1}}^0$ fail to form a generic arrangement. The arrangement $\B(n, k, \A)$ of hyperplanes $D_L$ is called $discriminantal$ $arrangement$ and has been introduced by Manin and Schechtman in \cite{MS}\footnote{Notice that Manin and Schechtman defined the discriminantal arrangement starting from a generic arrangement instead of its central translated as we do in this paper. For our purpose the latter is a more convenient choice.}. Notice that $\B(n, k, \A)$ depends on the trace at infinity $\A_{\infty}$ hence it is sometimes more properly denoted by $\B(n, k,\mathcal{A}_{\infty})$.

\subsection{Very generic and non-very generic discriminantal arrangements}\label{subs:lattice}
 It is well known (see, among others \cite{Crapo},\cite{MS}) that there exists an open Zarisky set $\mathcal{Z}$ in the space of (central) generic arrangements of $n$ hyperplanes in $\CC^k$, such that the intersection lattice of the discriminantal arrangement $\mathcal{B}(n,k,\A)$ is independent from the choice of the arrangement $\A \in  \mathcal{Z}$. Bayer and Brandt in \cite{BB} call the arrangements $\A \in  \mathcal{Z}$ \textit{very generic} and the ones which are not in $\mathcal{Z}$, \textit{non-very generic}. We will use their terminology in the rest of this paper. The name very generic comes from the fact that in this case the cardinality of the intersection lattice of $\B(n,k,\A)$ is the largest possible for any (central) generic arrangement $\A$ of $n$ hyperplanes in $\CC^k$. \\
In \cite{Crapo} Crapo proved that the intersection lattice of $\mathcal{B}(n,k,\A), \A \in  \mathcal{Z}$ is isomorphic to the Dilworth completion of the $k$-times lower-truncated Boolean algebra $B_n$ (see Theorem 2. page 149). A more precise description of this lattice is due to Athanasiadis who proved in \cite{Atha} a conjecture by Bayer and Brandt which stated that the intersection lattice of the discriminantal arrangement in very generic case is isomorphic to the collection of all sets $\{S_1, \ldots, S_m\}$, $S_i$ $\subset$ $[n]=\{1,\ldots,n\}$, $\mid S_i \mid \geq k+1$, such that
\begin{equation}\label{eq:vgcon}
\mid \bigcup_{i \in I} S_i \mid > k + \sum_{i \in I}(\mid S_i \mid - k) \mbox{ for all } I \subset [m]=\{1,\ldots,m\}, \mid I \mid \geq 2 \quad .
\end{equation}
 The isomorphism is the natural one which associate to the set $S_i$ the space $D_{S_i}=\bigcap_{L \subset S_i, \mid L \mid=k+1} D_L, D_L \in \mathcal{B}(n,k,\A)$ of all translated  of $\A$ having hyperplanes indexed in $S_i$ intersecting in a not empty space. In particular $\{S_1, \ldots, S_m\}$ will correspond to the intersection $\bigcap_{i=1}^m D_{S_i}$. \\
If $\A$ is very generic and the condition (\ref{eq:vgcon}) is satisfied, then the subspaces $D_{S_i}, i=1,\ldots, m$ intersect transversally (Corollary 3.6 in \cite{Atha}) or, equivalently, since $\rank~D_{S_i} = \mid S_i \mid - k$, that 
\begin{equation}\label{eq:Crapocond}
\rank \bigcap_{i=1}^m D_{S_i}=\sum_{i=1}^m ( \mid S_i \mid -k )
\end{equation}
Notice that, if $\A$ is very generic, the condition (\ref{eq:vgcon}) is always satisfied (see also \cite{Atha}) if
\begin{equation}\label{eq:crapo1} 
\bigcap_{i \in I}D_{S_i} \neq D_S, \mid S \mid > k+1 \mbox{ for any } I \subset [r]=\{1,\ldots, m\}, \mid I \mid \geq 2. 
\end{equation}
The fact that the condition (\ref{eq:crapo1}) implies the condition (\ref{eq:Crapocond}) for any set $\{S_1, \ldots, S_m\}$ which satisfies the condition (\ref{eq:vgcon}) corresponds to the definition provided by Crapo in \cite{Crapo} of a geometry of circuits\footnote{Here Crapo followed the preference of his advisor Rota who rarely used the name matroid.}. 


\noindent
Contrary to the very generic case, very few is known about the non-very generic case. In recent papers the first author (see \cite{LS}) and the first and second authors (see \cite{SSY1}, \cite{SSY2}) showed that the non-very generic arrangements are arrangements which hyperplanes give rise to special configurations (e.g. Pappus's configuration or Hesse configuration). Following this direction, in the rest of the paper we further develop the result in \cite{LS} providing a geometric sufficient condition for a central generic arrangement $\A$ to be non-very generic.
In order to do this we provide the following definition.

\begin{defi}\label{def:simple}An element $X$ in the intersection lattice of the discriminantal arrangement $\mathcal{B}(n,k,\A)$ is said to be a \textbf{simple} intersection if $X=\bigcap_{i=1}^r D_{L_i}, |L_i| = k+1$ and $\bigcap_{i \in I}D_{L_i} \neq D_S, \mid S \mid > k+1$ for any $I \subset [r], \mid I \mid \geq 2$. We call the number $r$ of the hyperplanes intersecting in $X$ the \textbf{multiplicity} of the simple intersection $X$ . 
\end{defi}
\noindent
The above considerations lead to the following Proposition.
\begin{prop}\label{pro:main}If the intersection lattice of the discriminantal arrangement $\mathcal{B}(n,k,\A)$ contains a simple intersection of rank strictly less than its multiplicity, then $\A$ is non-very generic.
\end{prop}

\proof Let's assume, by absurd, that $\A$ is a very generic arrangement such that $\mathcal{B}(n,k,\A)$ contains a simple intersection $X=\bigcap_{i=1}^r D_{L_i}, |L_i| = k+1$ of rank $s<r$. Then by definition of simple intersection, the set $\{L_1,\ldots, L_r \}$ satisfies the equation (\ref{eq:vgcon}) which, in turns, implies that the equation (\ref{eq:Crapocond}) is satisfied, i.e. rank $X=\sum_{i=1}^r ( \mid L_i \mid -k )=r$, which is an absurd and the proof is concluded. \endproof 

\noindent
Proposition \ref{pro:main} will play an important role in the rest of the paper since we will focus on a necessary and sufficient condition for the existence of such a simple intersection.

\subsection{Motivating examples}\label{ex:moti}
In this subsection we provide the two main examples given by Crapo (see \cite{Crapo}) and Falk (see \cite{CR} and \cite{Falk}) of simple but non-very generic intersections. Those two examples inspired the rest of the content of this paper.

\noindent
\textbf{Crapo example} is illustrated in Figure \ref{fig:quad}. In this case $\A^0$ is the central generic arrangement in the top of the figure while the arrangement $\A^{t}$ on the right of Figure \ref{fig:quad} is an element in the simple intersection $X=\bigcap_{i=1}^4 D_{L_i}$ with $L_1 = \{ 1,2,3 \}$, $L_2 = \{ 1,4,5 \}$, $L_3 = \{ 2,4,6 \}$, $L_4 = \{ 3,5,6 \}$. By definition of discriminantal arrangement, the only rank 4 intersection of $\B(6,2, \A^0)$ is given by $D_{[6]}$, the space of all central translated of $\A^0$. Since $\A^t$ is not central, this implies that $ X \neq D_{[6]}$ and hence $\rank~X< 4$, that is $X$ is a simple intersection of multiplicity 4 and $\rank~3 < 4$. That is the arrangement $\A^0$ in Figure \ref{fig:quad} is non-very generic.

\begin{center}
\begin{figure}[h]
\includegraphics[width=10cm]{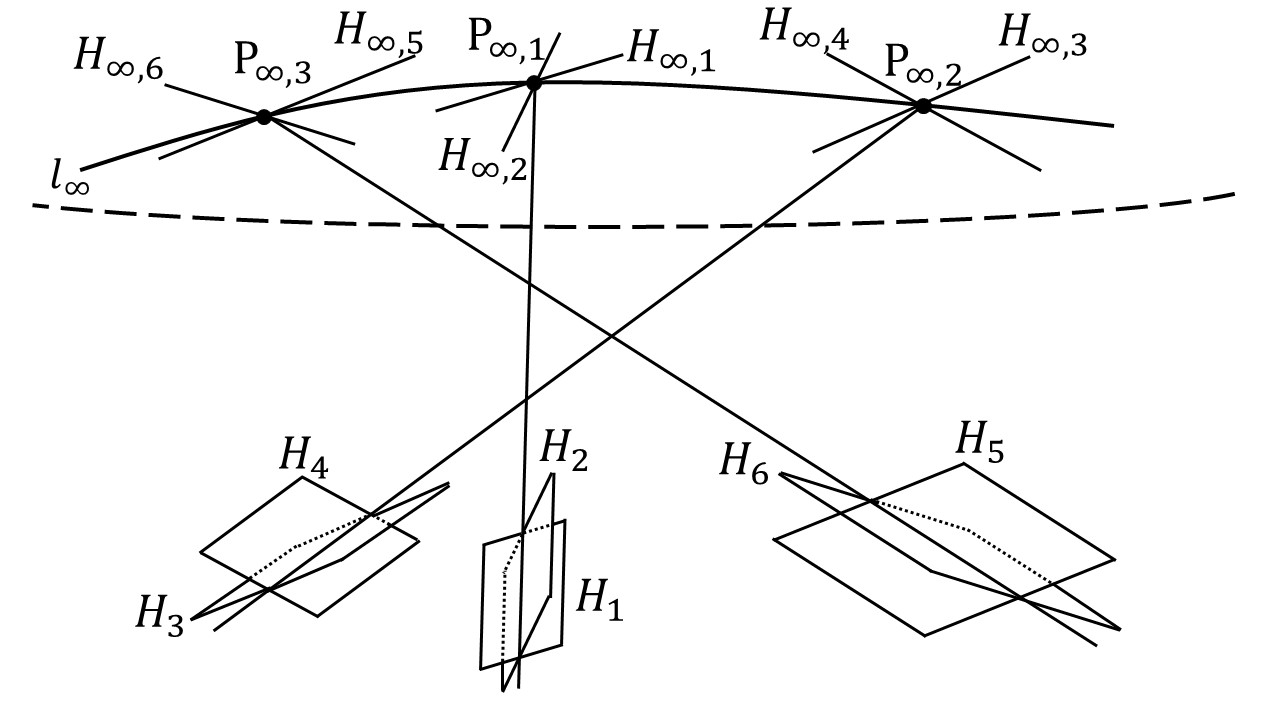}
\caption{Figure of non-very generic arrangement with 6 hyperplanes in $\RR^3$.}\label{fig:falk}
\end{figure}
\end{center}

\noindent
\textbf{Falk example} is illustrated in Figure \ref{fig:falk}. Let $\B(6,3,\A_\infty)$ be the discriminantal arrangement associated to a central generic arrangement $\A^0$ of 6 hyperplanes $H_i$ in $\RR^3$ which satisfy the condition that $\overline{H_i} \cap \overline{H_{i+1}} \cap H_\infty$, $i = 1,3,5$ span a line at infinity (see Figure \ref{fig:falk}). In \cite{LS} the authors proved that such an arrangement $\A^0$ admits a translation $\A^t$ which belongs to the simple intersection $X=\bigcap_{i=1}^3 D_{L_i}$ with $L_1 = \{ 1,2,3,4 \}$, $L_2 = \{ 1,2,5,6 \}$, $L_3 = \{ 3,4,5,6 \}$, that is, in particular, $\A^t$ is not a central arrangement. By definition of discriminantal arrangement, the only element in rank 3 in $\B(6,3, \A^0)$ is $D_{[6]}$ the space of all central translated of $\A^0$ hence $\rank~X< 3$, that is $X$ is a simple intersection of multiplicity 3 and $\rank~2 < 3$, i.e. $\A^0$ is non-very generic.
\section{A geometric condition for non-very genericity}\label{sec:geom}
In this section we provide a necessary and sufficient condition for the existence of a simple intersection $X$ of multiplicity $r$ and $\rank~X<r$ in the intersection lattice of the discriminantal arrangement $\mathcal{B}(n,k,\A^0)$. By Proposition \ref{pro:main} this is a sufficient condition for $\A^0$ to be non-very generic.

\begin{nota}\label{nota:fix}
To begin with let us fix some notations we will use throughout the paper.
\begin{itemize}
\item $\A^0$ is a central generic arrangement of $n$ hyperplanes in $\CC^k$\\
\item For each subset $L$ of $\{ 1, \dots, n \}$ with $|L| = k+1$, $D_L \subset \CC^n$ will denote the hyperplane in $\B(n,k, \A^0)$ corresponding to the subset $L$. \\
\item Fixed a set $\TT = \{ L_1, \dots, L_r \}$ of subsets $L_i \subset [n]$, $|L_i| = k+1$, for any arrangement $\A=\{H_1,\ldots, H_n\}$ translated of $\A^0$ we will denote by $P_i = \bigcap_{p \in L_i} H_p$ and $H_{i,j} = \bigcap_{p \in L_i \cap L_j} H_p$. Notice that $P_i$ is a point if and only if $\A \in D_{L_i}$, it is empty otherwise.\\
\end{itemize}
\end{nota}    
\subsection{$\mathbf{K_{\TT}}$-translated and $\mathbf{K_{\TT}}$-configurations}\label{sub:KT}
Let $\A^0=\{H_1^0, \ldots, H_n^0\}$ be a central generic arrangement in $\CC^k$, $\TT = \{ L_1, \dots, L_r \}$ fixed as in Notation \ref{nota:fix} and such that the conditions
 \begin{equation}\label{eq:proper1}
 \bigcup_{i=1}^r L_i = \bigcup_{i \in I} L_i  \quad \mbox{ and} \quad L_i \cap L_j \neq \emptyset
 \end{equation}
 are satisfied for any subset $I \subset [r], \mid I \mid=r-1$ and any two indices $1 \leq i < j \leq r$. In the rest of the paper a set $\TT$ which satisfies those properties will be called an $r$-\textbf{set}.\\
Given an $r$-set $\TT$ we provide the following two definitions.\\
\paragraph{$\mathbf{K_{\TT}}$-\textbf{translated}} A translated $\A=\{H_1, \ldots, H_n\}$ of $\A^0$ will be called $\mathbf{K_{\TT}}$ or $\mathbf{K_{\TT}}$-translated if $$\bigcap_{p \in L_i} H_p \neq \emptyset \mbox{ and } \bigcap_{p \in L_i \cup \{t\}} H_p = \emptyset$$ for any $i \in [r]$ and $t \notin L_i$.
\paragraph{$\mathbf{K_{\TT}}$-\textbf{configuration} $K_\TT(\A)$}The complete graph (as depicted in Figure \ref{fig:K_T_conf}) having the points $P_i=\bigcap_{p \in L_i} H_p$ as vertices and  the vectors $P_iP_j$ joining $P_i$ and $P_j$ as edges will be called \textbf{$K_{\TT}$-configuration} and denoted by $K_\TT(\A)$ (examples of graphs $K_\TT(\A)$ for  $|\TT| = 3,4,5$ are represented in Figure \ref{fig:K_T_conf_ex}). Notice that $P_iP_j \in H_{i,j} = \bigcap_{p \in L_i \cap L_j} H_p \neq \emptyset$ for any $1 \leq i<j\leq r$.\\

\noindent
For our purpose it is convenient to introduce a slightly weaker notion, beside the one of $\mathbf{K_{\TT}}$-translated and $\mathbf{K_{\TT}}$-configuration as follows. 
\begin{itemize}   
\item $\A$ will be called \textbf{almost} $\mathbf{K_{\TT}}$ if it is $\mathbf{K_{\TT}}$ but for one hyperplane $H_l^0$ and a set $S_l$, i.e. if there exists an hyperplane $H_l \in \A, l \in \bigcup_{i=1}^r L_i \setminus \bigcap_{i=1}^r L_i$, and a set $S_l \subseteq \{ L_i \in \TT \mid l \in L_i \}$ such that $\bigcap_{p \in L_j} H_p \neq \emptyset$ for any $L_j \in \TT \setminus S_l$ and $\bigcap_{p \in L_i} H_p = \emptyset$ for any $L_i \in S_l$.
 \item If we keep the notation $P_i=\bigcap_{p \in L_i \setminus \{ l \}} H_p$, $L_i \in S_l, P_i = \bigcap_{p \in L_i} H_p, L_i \notin S_l$, the complete graph having $P_i$ as vertices and $P_iP_j$ as edges will be called \textbf{almost $K_{\TT}$-configuration} and denoted by $K_{\TT \setminus S_l} (\A)$.\\
\end{itemize}
Notice that since $\A^0$ is a central generic arrangement in $\C^k$ and $\mid L_i \mid = k+1$, then $\bigcap_{p \in L_i \setminus \{ l \}} H_p \neq \emptyset$ for any $L_i \in S_l$. Moreover if the set $\{ L_i \in \TT \mid l \in L_i \}$ is not empty then its cardinality is $\mid \{ L_i \in \TT \mid l \in L_i \} \mid \geq 2$ since $l$ must belong to at least two sets inside $\TT$ by the first condition in equation (\ref{eq:proper1}) . 

\begin{center}
\begin{figure}[h]
\begin{tikzpicture}[line join = round, line cap = round]
\coordinate [label=above:$P_1$] (0) at (0,0);
\coordinate [label=above:$P_2$] (1) at ({-sqrt(3)*(3/5)},{-1*(3/5)});
\coordinate [label=below:$P_3$] (2) at ({-sqrt(3)*(3/5)},{-3*(3/5)});
\coordinate [label=below:$\dots$] (3) at (0,{-4*(3/5)});
\coordinate [label=below:$\dots$] (4) at ({sqrt(3)*(3/5)},{-3*(3/5)});
\coordinate [label=above:$P_r$] (5) at ({sqrt(3)*(3/5)},{-1*(3/5)});
\coordinate [label=\reflectbox{$\ddots$}] (6) at ({sqrt(3)/2 *(3/5)},{-4*(3/5)});
\begin{scope}
\draw (0) -- node[above] {$H_{1,2}$} (1);
\draw (0) --  (2);
\draw (0) --  (4);
\draw (0) -- node[above] {$H_{1,r}$} (5);
\draw (1) -- node[left] {$H_{2,3}$} (2);
\draw (1) --  (3);
\draw (1) --  (5);
\draw (2) --  (3);
\draw (2) --  (5);
\draw (4) --  (5);
\fill (0) circle (2pt);
\fill (1) circle (2pt);
\fill (2) circle (2pt);
\fill (5) circle (2pt);
\end{scope}
\end{tikzpicture}
\caption{$K_\TT$-configuration for $|\TT| = r$}\label{fig:K_T_conf}
\end{figure}
\end{center}

\begin{center}
\begin{figure}[h]
\begin{tikzpicture}[line join = round, line cap = round]
\coordinate [label=below:$P_1$] (0) at (0,0);
\coordinate [label=below:$P_2$] (1) at (4,0);
\coordinate [label=above:$P_3$] (2) at (2, {2*sqrt(3)});

\begin{scope}
\draw (0)-- node[below] {$H_{1,2}$} (1);
\draw (0)-- node[left] {$H_{1,3}$} (2);
\draw (1)-- node[right] {$H_{2,3}$} (2);
\fill (0) circle (2pt);
\fill (1) circle (2pt);
\fill (2) circle (2pt);
\end{scope}
\end{tikzpicture}
\begin{tikzpicture}[line join = round, line cap = round]
\coordinate [label=below:$P_1$] (0) at (0,0);
\coordinate [label=below:$P_2$] (1) at (4,0);
\coordinate [label=above:$P_3$] (2) at (4,4);
\coordinate [label=above:$P_4$] (3) at (0,4);

\begin{scope}
\draw (0) -- node[below] {$H_{1,2}$} (1);
\draw (0) --  (2);
\node[] at (1,1.5) {$H_{1,3}$};
\draw (0) -- node[left] {$H_{1,4}$} (3);
\draw (1) -- node[right] {$H_{2,3}$} (2);
\draw (1) --  (3);
\node[] at (3,1.5) {$H_{2,4}$};
\draw (2) -- node[above] {$H_{3,4}$} (3);
\fill (0) circle (2pt);
\fill (1) circle (2pt);
\fill (2) circle (2pt);
\fill (3) circle (2pt);
\end{scope}
\end{tikzpicture}
\begin{tikzpicture}
\coordinate [label=below:$P_1$] (0) at (0,0);
\coordinate [label=below:$P_2$] (1) at (2,0);
\coordinate [label=right:$P_3$] (2) at (3, 5/2);
\coordinate [label=above:$P_4$] (3) at (1, 4);
\coordinate [label=left:$P_5$] (4) at (-1, 5/2);

\begin{scope}
\draw (0)-- node[below] {$H_{1,2}$} (1);
\draw (0)--  (2);
\node[] at (2.1,1.5) {$H_{1,3}$};
\draw (0)--  (3);
\node[] at (0.4,2) {$H_{1,4}$};
\draw (0)-- node[left] {$H_{1,5}$} (4);
\draw (1)-- node[right] {$H_{2,3}$} (2);
\draw (1)--  (3);
\node[] at (1.6,2) {$H_{2,4}$};
\draw (1)--  (4);
\node[] at (0,1.5) {$H_{2,5}$};
\draw (2)-- node[above] {$H_{3,4}$} (3);
\draw (2)--  (4);
\node[] at (1.8,2.7) {$H_{3,5}$};
\draw (3)-- node[above] {$H_{4,5}$} (4);
\fill (0) circle (2pt);
\fill (1) circle (2pt);
\fill (2) circle (2pt);
\fill (3) circle (2pt);
\fill (4) circle (2pt);
\end{scope}
\end{tikzpicture}
\caption{Left: $|\TT|=3$, Center: $|\TT|=4$, Right: $|\TT|=5$}\label{fig:K_T_conf_ex}
\end{figure}
\end{center}

\begin{rem} Similarly to the $K_{\TT}$-configuration we could define the $\Delta_{\TT}$-configuration as the simplicial complex having as $t$-face $P_{i_1}\ldots P_{i_{t+1}} \in  \bigcap_{p \in \cap_{j=1}^{t+1} L_{i_j}} H_p \neq \emptyset$. 
Notice that  in general, the intersection $\bigcap_{p \in \cap_{j=1}^{t+1} L_{i_j}} H_p$ can be empty, that is $\Delta_{\TT}$ is not a simplex. As pointed out by Crapo in \cite{Crapo}, this simplicial complex may play a fundamental role in the study of non-very generic arrangements.
\end{rem}

\noindent
In order to understand the ratio behind the definition of $\mathbf{K_{\TT}}$ and \textbf{almost} $\mathbf{K_{\TT}}$-translated let's consider the Crapo's configuration depicted in Figure  \ref{fig:quad}.

\begin{figure}[h]
\centering
\begin{tikzpicture}
\coordinate (0) at (-1,0);
\coordinate (1) at (-1/2,-2/3);
\coordinate (2) at (-1/3,-1);
\coordinate (3) at (17/4,-1);
\coordinate (4) at (5,-1/2);
\coordinate (5) at (5,0);
\coordinate (6) at (4,9/2);
\coordinate (7) at (4,16/3);
\coordinate (8) at (11/4,5);
\coordinate (9) at (4/3,4);
\coordinate (10) at (1/2,4);
\coordinate (12) at (0,5/2);
\coordinate [label=below:$P_1$] (p1) at (0,0);
\coordinate [label=below:$P_2$] (p2) at (4,0);
\coordinate [label=$P_3$] (p3) at (3,4);
\coordinate [label=$P_4$] (p4) at (1,3);


\begin{scope}
\draw (p1) -- (p2);
\draw (p1) -- (p3);
\draw (p1) -- (p4);
\draw (p2) -- (p3);
\draw (p3) -- (p4);
\fill (p1) circle (1.5pt);
\fill (p2) circle (1.5pt);
\fill (p3) circle (1.5pt);
\fill (p4) circle (1.5pt);

\end{scope}
\end{tikzpicture} \ \ \ \ \ 
\begin{tikzpicture}
\coordinate (0) at (-1,0);
\coordinate (1) at (-1/2,-2/3);
\coordinate (2) at (-1/3,-1);
\coordinate (3) at (17/4,-1);
\coordinate (4) at (9/2,-1/2);
\coordinate (5) at (5,0);
\coordinate (6) at (4,9/2);
\coordinate (7) at (4,16/3);
\coordinate (8) at (11/4,5);
\coordinate (9) at (4/3,4);
\coordinate (10) at (0,4);
\coordinate (12) at (0,5/2);
\coordinate [label=below:$P_1$] (p1) at (0,0);
\coordinate [label=below:$P_2$] (p2) at (4,0);
\coordinate [label=$P_3$] (p3) at (3,4);
\coordinate [label=$P_4$] (p4) at (1,3);

\begin{scope}
\draw (p1) -- (p2);
\draw (p1) -- (p3);
\draw (p1) -- (p4);
\draw (p2) -- (p3);
\draw (p2) -- (p4);
\draw (p3) -- (p4);
\fill (p1) circle (1.5pt);
\fill (p2) circle (1.5pt);
\fill (p3) circle (1.5pt);
\fill (p4) circle (1.5pt);

\end{scope}
\end{tikzpicture}
\caption{$K_{\TT \setminus S_5}(\A)$, $S_5 = \{ L_2, L_4\}$ on the left and $K_\TT(\A')$ on the right.}\label{pic:Crapo_ex_gr}
\end{figure}
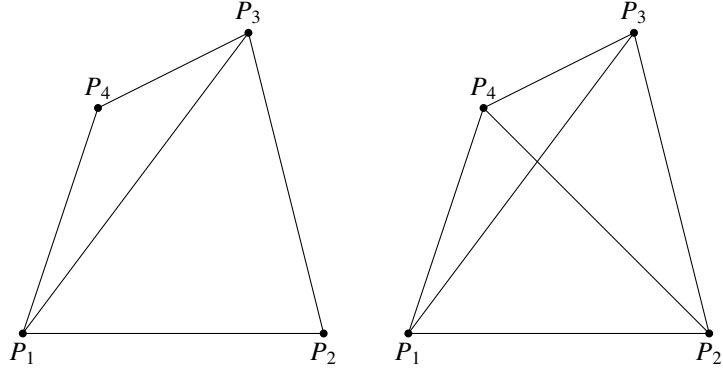

\begin{ex}[Crapo's example]\label{ex:crapo2}
Let us consider the Crapo's example in Subsection \ref{ex:moti}. The non-very genericity of Crapo's arrangement $\A^0$ implies that any translation $\A$ of $\A^0$ such that $\A \in \bigcap_{i=1}^3 D_{L_i}$ has to satisfy $\A \in D_{L_4}$. 
That is, as depicted in Figure \ref{pic:Crapo_ex}, any translation $t_5$ of the hyperplane $H_5^0$ for which $P_2 \in H_5^{t_5}$ has to satisfy $P_4 \in H_5^{t_5}$.\\
In other words, if we choose $l = 5$, then the \textbf{almost} $K_\TT$-configuration given by $P_1 = H^{t_1}_1 \cap H^{t_2}_2 \cap H_3^{t_3}$, $P_2 = H_1^{t_1} \cap H_4^{t_4}$, $P_3 = H_2^{t_2} \cap H_4^{t_4} \cap H_6^{t_6}$ and $P_4 = H_3^{t_3} \cap H_6^{t_6}$ becomes a $K_\TT$-configuration (see Figure \ref{pic:Crapo_ex_gr}) by the only translation $t_5$ of $H_5^0$ which satisfies $P_2 \in H_5^{t_5}$. 
\end{ex}

\noindent
The Example \ref{ex:crapo2} also motivated the following main definition.

\begin{defi}\label{df:dep} A central generic arrangement $\A^0$ of $n$ hyperplanes in $\CC^k$ is called \textbf{(r,s)-dependent} if there exist an $r$-set $\TT = \{ L_1, \dots, L_r \}$, an index $l \in \bigcup_{i=1}^r L_i \setminus \bigcap_{i=1}^r L_i$ and a subset $S_l \subseteq \{ L_i \in \TT \mid l \in L_i \}$, $\mid S_l \mid = s$, such that any \textbf{almost} $K_\TT$-configuration $K_{\TT \setminus S_l} (\A)$ gives rise to a $K_\TT$-configuration $K_\TT(\A')$ with the $K_\TT$-translated $\A'$ obtained from the \textbf{almost} $K_\TT$-translated $\A$  by a suitable translation of the hyperplane $H_l \in \A$. If $s = 2$, then we call $\A^0$ $r$-dependent.
\end{defi}

\begin{figure}[h]
\centering
\begin{tikzpicture}
\coordinate (0) at (-1,0);
\coordinate (1) at (-1/2,-2/3);
\coordinate (2) at (-1/3,-1);
\coordinate (3) at (17/4,-1);
\coordinate (4) at (5,-1/2);
\coordinate (5) at (5,0);
\coordinate (6) at (4,9/2);
\coordinate (7) at (4,16/3);
\coordinate (8) at (11/4,5);
\coordinate (9) at (4/3,4);
\coordinate (10) at (1/2,4);
\coordinate (12) at (0,5/2);
\coordinate (p1) at (0,0);
\coordinate (p2) at (4,0);
\coordinate (p3) at (3,4);
\coordinate (p4) at (1,3);

\coordinate [label=$H_1^{t_1}$] (H1) at (-1.2,-0.3);
\coordinate [label=$H_2^{t_2}$] (H2) at (-0.6,-1.1);
\coordinate [label=$H_3^{t_3}$] (H3) at (-0.3,-1.5);
\coordinate [label=$H_4^{t_4}$] (H4) at (4.3,-1.5);
\coordinate [label=$H_5^0$] (H5) at (5.2,-0.9);
\coordinate [label=$H_6^{t_6}$] (H6) at (4.3,4.3);
\coordinate [label=$P_1$] (P1) at (-0.2,-0.1);
\coordinate [label=$P_2$] (P2) at (4.2,-0.1);
\coordinate [label=$P_3$] (P3) at (2.7,3.9);
\coordinate [label=$P_4$] (P4) at (0.9,3.1);

\begin{scope}
\draw (0) -- (5);
\draw (1) -- (7);
\draw (3) -- (8);
\draw (2) -- (9);
\draw[dashed] (4) -- (10);
\draw (6) -- (12);
\fill (p1) circle (1.5pt);
\fill (p2) circle (1.5pt);
\fill (p3) circle (1.5pt);
\fill (p4) circle (1.5pt);

\end{scope}
\end{tikzpicture} \ \ \ \ \ 
\begin{tikzpicture}
\coordinate (0) at (-1,0);
\coordinate (1) at (-1/2,-2/3);
\coordinate (2) at (-1/3,-1);
\coordinate (3) at (17/4,-1);
\coordinate (4) at (9/2,-1/2);
\coordinate (5) at (5,0);
\coordinate (6) at (4,9/2);
\coordinate (7) at (4,16/3);
\coordinate (8) at (11/4,5);
\coordinate (9) at (4/3,4);
\coordinate (10) at (0,4);
\coordinate (12) at (0,5/2);
\coordinate (p1) at (0,0);
\coordinate (p2) at (4,0);
\coordinate (p3) at (3,4);
\coordinate (p4) at (1,3);

\coordinate [label=$H_1^{t_1}$] (H1) at (-1.2,-0.3);
\coordinate [label=$H_2^{t_2}$] (H2) at (-0.6,-1.1);
\coordinate [label=$H_3^{t_3}$] (H3) at (-0.3,-1.5);
\coordinate [label=$H_4^{t_4}$] (H4) at (4.3,-1.5);
\coordinate [label=$H_5^{t_5}$] (H5) at (4.7,-0.9);
\coordinate [label=$H_6^{t_6}$] (H6) at (4.3,4.3);
\coordinate [label=$P_1$] (P1) at (-0.2,-0.1);
\coordinate [label=$P_2$] (P2) at (4.2,-0.1);
\coordinate [label=$P_3$] (P3) at (2.7,3.9);
\coordinate [label=$P_4$] (P4) at (0.9,3.1);

\begin{scope}
\draw (0) -- (5);
\draw (1) -- (7);
\draw (3) -- (8);
\draw (2) -- (9);
\draw (4) -- (10);
\draw (6) -- (12);
\fill (p1) circle (1.5pt);
\fill (p2) circle (1.5pt);
\fill (p3) circle (1.5pt);
\fill (p4) circle (1.5pt);

\end{scope}
\end{tikzpicture}
\caption{\textbf{Almost} $K_\TT$-translated $\A$ on the left and $K_\TT$-translated $\A'$ on the right.}\label{pic:Crapo_ex}
\end{figure}
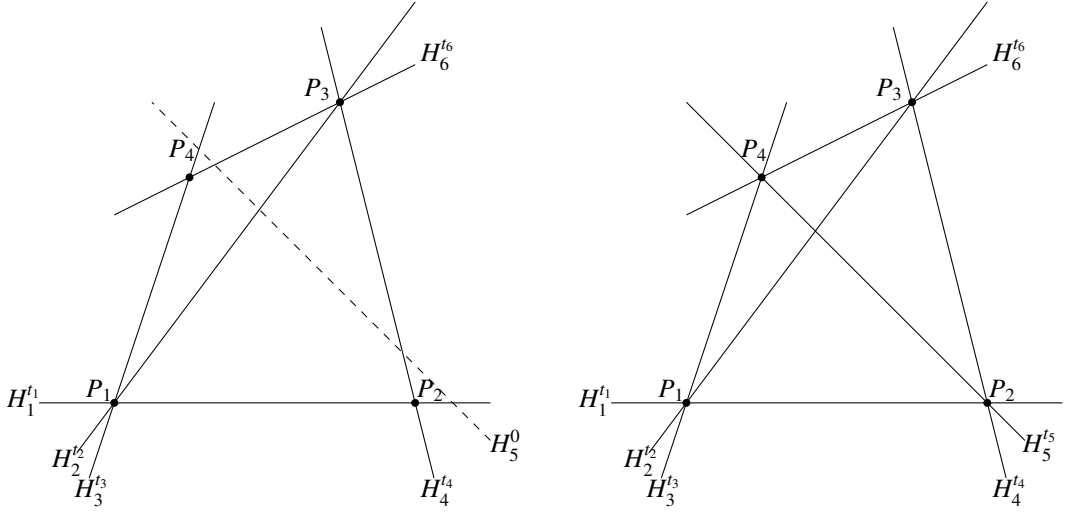

\noindent
By the argument in the Example \ref{ex:crapo2}, the Crapo's arrangement is $4$-dependent.\\
The definition of $(r,s)$-dependency generalizes the definition of dependency given in \cite{LS}. Indeed we have the following proposition which, in particular, applies to Falk's example (see Example \ref{ex2:fl}).

\begin{prop}A central generic arrangement $\A^0$ of $n$ hyperplanes in $\C^3$ is dependent if and only if it is $3$-dependent.
 \end{prop}
\proof
Let's consider the set $\TT = \{ L_1, L_2, L_3 \}$ and assume $n=3s, k=2s-1, s \geq 2$. The general case of any $n$ and $k$ is obtained considering a deletion and a restriction of the original arrangement as in \cite{LS}. \\
In order to prove the statement it is enough to show that both conditions, i.e. 3-dependency and dependency, are equivalent to the condition that the space $H_{i,j}$ is a subspace of $H_{i,k} + H_{k,j}$.

\textbf{Dependency.} Recall that an arrangement $\A^0$ of $3s$ hyperplanes in $\CC^{2s-1}$, $s \geq 2$ is dependent if it exists a set $\TT = \{ L_1, L_2, L_3 \}$ of subsets $L_i \subset [3s]$ such that $\mid L_i \mid = 2s$, $\mid L_i \cap L_j \mid = s$, $\mid \bigcup_{i=1}^3 L_i \mid = 3s$ and the spaces $H_{i,j}=\bigcap_{p \in L_i \cap L_j} H_{p}$ span a subspace of dimension $2s-2$ in $\CC^{2s-1}$. The condition $\mid L_i \cap L_j \mid = s$ implies that $H_{i,j}$ are spaces of codimension $s$, that is of dimension $s-1$ in $\CC^{2s-1}$. Moreover $\mid \bigcup_{i=1}^3 L_i \mid = 3s$ implies that the $\bigcup_{i=1}^3 L_i $ is disjoint union of the three sets $L_i \cap L_j$, that is any two subspaces $H_{i,j}$ are in direct sum, i.e. $H_{i,k} \oplus H_{k,j}$ span a space of dimension $2s-2$. Hence dependency condition is equivalent to the fact that $H_{i,j}$ belongs to the space generated by $H_{i,k} \oplus H_{k,j}$. 

\textbf{$3$-dependency.}  First of all notice that, when $n=3s, k=2s-1, s \geq 2$ the condition for which $\TT$ is a $3$-set implies that the subsets $L_i \subset [3s]$ satisfy $\mid L_i \mid = 2s$, $\mid L_i \cap L_j \mid = s$, $\mid \bigcup_{i=1}^3 L_i \mid = 3s$ and, in particular $\bigcap_{i=1}^3 L_i = \emptyset$. Moreover a $K_\TT$-configuration when $\mid \TT \mid = 3$ is equivalent to the fact that $P_iP_j=P_iP_k+P_kP_j$ (see Figure \ref{fig:K_T_conf_ex}) and the condition $\bigcup_{i=1}^3 L_i = \bigcup_{i \in I \subset [3], \mid I \mid = 2} L_i$ implies that any index $l \in \bigcup_{i=1}^3 L_i \setminus \bigcap_{i=1}^3 L_i$ belongs to exactly two different subsets $L_i$ and $L_j$. The $3$-dependency condition is then equivalent to the fact that any translation for which the vertex $P_i=\bigcap_{p \in L_i}H_p \neq \emptyset$ exists the vertex $P_j=\bigcap_{p \in L_j}H_p \neq \emptyset$ has to exist. Hence, in particular, $P_iP_j=P_iP_k + P_kP_j$ for any $P_i,P_j \in H_{i,j}$, that is $H_{i,j}$ is a subspace of $H_{i,k} \oplus H_{k,j}$ which are in direct sum since $\bigcap_{i=1}^3 L_i = \emptyset$.
\endproof

\begin{center}
\begin{figure}[h]
\includegraphics[width=10cm]{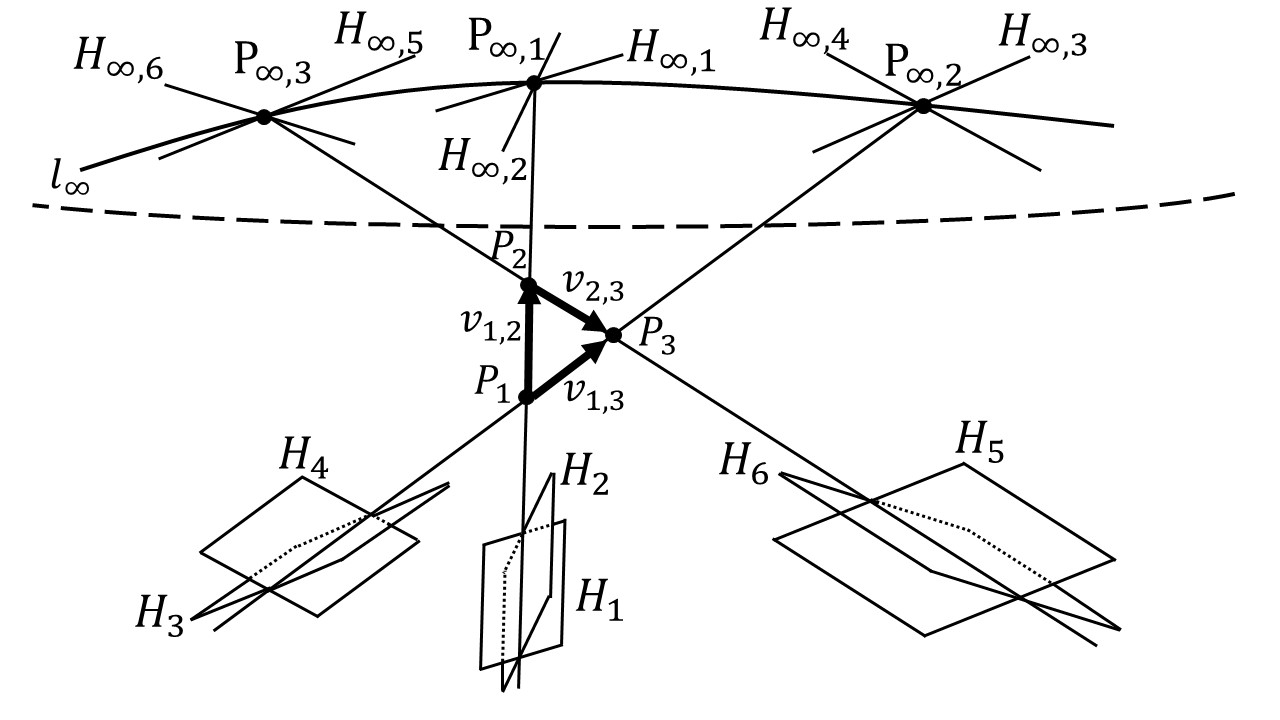}
\caption{}\label{fig:falk_3dep}
\end{figure}
\end{center}

\begin{ex}[Falk's example]\label{ex2:fl}
Consider the Falk's example in Subection \ref{ex:moti}. In this case $\A^0$ is an arrangement of 6 hyperplanes in $\RR^3$ and the set $\TT = \{ L_1, L_2, L_3 \}$ is given by $L_1 = \{ 1,2,3,4 \}$, $L_2 = \{ 1,2,5,6 \}$ and $L_3 = \{ 3,4,5,6 \}$ which satisfy the conditions $\mid L_i \mid = 4$, $\mid L_i \cap L_j \mid = 2$ and $\mid \bigcup_{i=1}^3 L_i \mid = 6$. Since the spaces $H_{i,j, \infty}=\bigcap_{p \in L_i \cap L_j} \overline{H}_{p} \cap H_\infty$ span a subspace of dimension $1$ in $H_\infty$ then $H_{i,j}$ span a space of dimension $2$ in $\RR^3$, that is $\A^0$ is a dependent arrangement.\\
In other words, let's choose the index $l=6$ and $S_l=\{L_2,L_3\}$. Then for any translated $\A$ of $\A^0$ such that $P_1 = \bigcap_{p \in L_1}H_p \neq \emptyset$ exists, the translated $H'_6$ of $H_6$ for which $P_2 = H'_6 \cap \bigcap_{p \in L_2 \setminus \{ 6 \}}H_p \neq \emptyset$ also satisfies $P_3 = H'_6 \cap \bigcap_{p \in L_3 \setminus \{ 6 \}}H_p \neq \emptyset$. Moreover $P_1P_3 = P_1P_2 + P_2P_3$, that is $\A^0$ is $3$-dependent (see Figure \ref{fig:falk_3dep}).
\end{ex}

\noindent
Notice that the condition of $r$-dependency is non trivial. Indeed by $\bigcup_{i=1}^r L_i = \bigcup_{i \in I \subset [r], \mid I \mid=r-1} L_i$ it follows that any index $l \in \bigcup_{i=1}^r L_i$ has to belong to at least two different subsets $L_i$'s. Hence if $L_i \neq L_j$ are two different subsets containing the index $l$ the fact that $H_l$ is a translated of $H^0_l$ for which $\bigcap_{p \in L_i\setminus \{l\}} H_p \cap H_l \neq \emptyset$ does not imply, in general, that $\bigcap_{p \in L_j\setminus \{l\}} H_p \cap H_l \neq \emptyset$. In particular $(r,s)$-dependency is always non trivial for any $1 < s \leq \mid \{L_i \in \TT \mid l \in L_i\} \mid$. The following lemma holds.

\begin{lem}\label{lem:almostKT}
If a central generic arrangement $\A^0$ of $n$ hyperplanes in $\C^k$ is $(r,s)$-dependent for some $s > 1$ then the discriminantal arrangement $\B(n,k,\A^0)$ admits a simple intersection $X$ of multiplicity $r$ and $\rank~X<r$. 
\end{lem}

\proof By definition $\A^0$ is $(r,s)$-dependent if and only if there exist an $r$-set $\TT = \{ L_1, \dots, L_r \}$, an index $l \in \bigcup_{i=1}^r L_i \setminus \bigcap_{i=1}^r L_i$ and a subset $S_l \subseteq \{ L_i \in \TT \mid l \in L_i \}$, $\mid S_l \mid = s>1$, such that any translated arrangement $\A \in \bigcap_{L_i \in \TT \setminus S_l} D_{L_i}$ that satisfy $\A \in D_{L_j}$, for an $L_j \in S_l$, has to satisfy $\A \in \bigcap_{L_i \in \TT} D_{L_i}=X$. That is, if $\A^0$ is $(r,s)$-dependent then 
\begin{equation}\label{eq:simp_nvg}
X=\bigcap_{L_i \in \TT} D_{L_i}= \bigcap_{L_i \in \TT \setminus {S_l}} D_{L_i} \cap D_{L_j}  
 \end{equation}
is a simple intersection of multiplicity $r$ and $\rank~X \leq r-s+1<r$ since $s>1$.
\endproof

\noindent
An immediate consequence of the Lemma \ref{lem:almostKT} and the Proposition \ref{pro:main} is the following main theorem.

\begin{thm}\label{th:conf}If a central generic arrangement $\A^0$ of $n$ hyperplanes in $\CC^k$ is $(r,s)$-dependent for some $s > 1$ then $\A^0$ is non-very generic.
\end{thm}

\noindent
The authors conjecture that, if the sets $L_i$'s in the Definition \ref{def:simple} are replaced by subsets $S_i \in [n]$ of any cardinality, then the argument in this paper can be rewritten to obtain a necessary and sufficient condition for a central generic arrangement to be
non-very generic. But while it is not difficult to generalize the definition of $K_\TT$-translated, how to define the analogous of \textbf{almost} $K_\TT$ and hence the equivalent of the $(r,s)$-dependency requires farther studies. Moreover the combinatorial condition on the set of indices in equation (\ref{eq:proper1}) has to be generalized. \\
Finally, let's remark that the geometric condition in Theorem \ref{th:conf} can be translated in a condition on the vectors belonging to the subspaces $H_{i,j}$. This condition, in turns, lead to a computational feasible way to build non-very generic arrangements.  This is the content of a forthcoming paper.


\end{document}